\documentclass[12pt]{article}
\usepackage{epsfig}
\usepackage{epstopdf}
\usepackage{graphics}
\usepackage{amsmath}
\usepackage{amssymb}

\newtheorem{thm}{Theorem}[section]

\newtheorem{hyp}[thm]{Hypothesis}

\newcommand{\Sone}{{\bf S}^1}

\newcommand{\SEtwo}{{\bf SE}(2)}
\newcommand{\Ttwo}{{\bf T}^2}
\newcommand{\Tthree}{{\bf T}^3}

\newcommand{\qed}{\hfill\mbox{\raggedright\rule{.07in}{.1in}}
  \vspace{1ex}} 

\newcommand{\Section}[1]{\section{#1} \setcounter{equation}{0}}

\newcommand{\arraystart}{\renewcommand{\arraystretch}{1.6}}
\newcommand{\arrayfinish}{\renewcommand{\arraystretch}{1.2}}
 
\title{Lattice symmetry breaking perturbations for spiral waves}
 
\author{Laurent Charette\footnote{Current address Department of Mathematics, University of British Columbia, 121-1984 Mathematics Road,
Vancouver BC, V6K 1T5 CANADA} \\Department of Mathematics and Statistics\\
University of Ottawa\\Ottawa, ON K1N 6N5\\CANADA  
\and Victor G. LeBlanc\\Department of
  Mathematics and Statistics\\University of Ottawa\\Ottawa, 
ON K1N 6N5\\CANADA}
\date{\today}

\begin{document}

\maketitle

\begin{abstract}
Spiral waves in two-dimensional excitable media have been observed experimentally and studied extensively.  It is now well-known that the symmetry properties of the medium of propagation drives many of the dynamics and bifurcations which are experimentally observed for these waves.  Also, symmetry-breaking induced by boundaries, inhomogeneities and anisotropy have all been shown to lead to different dynamical regimes as to that which is predicted for mathematical models which assume infinite homogeneous and isotropic planar geometry.  Recent mathematical analyses incorporating the concept of forced symmetry-breaking from the Euclidean group of all planar translations and rotations have given model-independent descriptions of the effects of media imperfections on spiral wave dynamics.  In this paper, we continue this program by considering rotating waves in dynamical systems which are small perturbations of a Euclidean-equivariant dynamical system, but for which the perturbation preserves only the symmetry of a regular square lattice.
\end{abstract}

\pagebreak
\Section{Introduction}
Two dimensional spiral waves are observed in many different situations \cite{AMN99, Bar92, Bar94, BKT90, BLM, DPSBJ92, GLM97, GLM98, Grill96, JW, LeBlanc, LW, Mesin, MZ, Rabinov, sandetal2in96, sandetalin98, Yuanetal}. In excitable biological tissue such as the cerebral cortex, the myocardium or the retina, spirals are usually precursors to serious pathological conditions \cite{KeenerSneyd}.  Thus there is great interest in understanding the dynamics of spiral waves, and how these are affected by model parameters and interactions with medium imperfections.

One of the most important theoretical breakthroughs in the study of spiral waves has been the adoption of techniques from group-equivariant dynamical systems \cite{Bar94, GLM97, GLM98, sandetal2in96, sandetalin98}.  Typically the physical systems in which spirals are observed are modeled using systems of reaction-diffusion partial differential equations
\begin{equation}\label{rdpde1}
\frac{\partial u}{\partial t}=D\cdot \nabla^2 u + f(u,\alpha)
\end{equation}
where $u(x,y,t)$ is a function $u:\mathbb{R}^2\times\mathbb{R}^{+}\longrightarrow\mathbb{R}^n$ which can represent a vector of concentrations, electrical potentials, {\it et cetera} depending on the system being modeled.  The function $f$ models local reaction dynamics, $D$ is an $n\times n$ matrix of diffusion coefficients, and $\alpha$ are model parameters.  

Equation (\ref{rdpde1}) is invariant under the full Euclidean group of planar rotations and translations, $\SEtwo$ (in fact, they are also invariant under reflections in space, but for our purposes, this symmetry will not be useful).  Rigidly and uniformly rotating spiral waves are now well-understood to be {\it rotating wave} solutions of (\ref{rdpde1}) in the sense of equivariant dynamical systems \cite{GSS88}; that is, evolution in time is the same as a uniform spatial rotation of the initial condition about some point in space.  This fact has led to a model-independent understanding of many of the observed dynamical states and bifurcations of spirals, e.g. transition to meandering, resonant linear drifting \cite{Bar94, GLM97, GLM98, sandetal2in96, sandetalin98}.

For the purposes of applications, equation (\ref{rdpde1}) and the associated symmetry describes a highly idealized situation which is almost never achieved in reality: infinite planar domain which is both homogeneous and isotropic.  In practice, media of propagation are finite in extent, so boundaries can play a significant role.  Also, especially in excitable biological tissue, there are many potential sources of anisotropy and inhomogeneities.  Many experiments (both numerical and physical) have shown that these imperfections have a significant impact on the dynamics of spiral waves \cite{Mesin, MZ91, MMV98, Roth1, Roth2, YP86, ZM96}

Using a model-independent dynamical systems formulation based on the idea of forced symmetry-breaking from $\SEtwo$ to a proper subgroup of $\SEtwo$, many of these effects (which result from inhomogeneities and/or anisotropy)
can be explained in terms of simple dynamics and bifurcations \cite{BLM, LW, LeBlanc, LeBlancRoth}.  In some cases, this approach has even led to some predictions which have been later confirmed experimentally \cite{BEL, LeBlancRoth}.

It is this forced symmetry-breaking approach that we wish to pursue in this paper.  Specifically, we will be interested in cases in which the Euclidean symmetry of (\ref{rdpde1}) is broken by a small perturbation which preserves only the symmetry of a regular square lattice.   At the cellular level,  excitable biological tissues exhibit inhomogeneities in conduction via the gap junctions which separate the cells \cite{KeenerSneyd}.  In cases where the spatial scale of the spiral wave is comparable to the cellular scale, it is not unreasonable to assume that the conduction inhomogeneities caused by these gap junctions may have an impact on the dynamics of the spiral.  Moreover, if one wants to investigate these effects, a reasonable simplifying first order approximation is to assume that the cells are distributed in a regular lattice array.  See also \cite{Adamatsky, Agladze, Sole, Xu} for further motivation.  Based on the results of \cite{BLM, LeBlanc, LW}, we expect that the lattice may have stabilising (anchoring) effects on the spirals and/or induce meandering.  These expectations will indeed be consequences of our results in this paper.

At some phenomenological level, it is also expected that our results may be a first step in understanding potential numerical effects induced by spatial-discretization of the domain for integrations of reaction-diffusion partial differential equations.  Of course in most situations, we expect that numerical discretization effects will be negligeable if the spatial grid is fine enough.  However, coarse grids (relative to the spatial scale of the spiral) could have effects on the observed dynamics, and our work here could point to what these effects may be.  In a specific instance, we argue that forced lattice symmetry-breaking gives results which are consistent with numerically-observed transition from rigidly rotating wave to linearly translating waves (with retracted tip) as the rotation frequency of the spiral wave tends to zero \cite{JW, MZ}.

This paper is organized as follows.  In the next section, we establish the symmetry properties we wish to consider, and derive the main center-bundle equations which represent the dynamics on perturbed relative equilibria (i.e. rotating or travelling waves) in the context of lattice forced symmetry-breaking.  In section 3, we recall a classic result by Hale on averaging and the existence of invariant manifolds for perturbed systems, in a form suitable for our purposes.  Section 4 deals with the effects of lattice symmetry-breaking perturbations on rotating spiral waves, and section 5 considers the effects of the perturbation on travelling waves.  We illustrate some of our results with numerical simulations which are presented in section 6.  The proofs of our two main results are presented in appendices A and B.

\Section{Preliminaries}

Let $\SEtwo$ denote the set of all planar translations and rotations, parametrized as follows
\[
\SEtwo=\{(\theta,p)\in\Sone\times\mathbb{R}^2\}
\]
and whose action on a point ${\displaystyle z=\left(\begin{array}{c}x\\y\end{array}\right)\in\mathbb{R}^2}$, is given by
\begin{equation}
(\theta,p)\cdot z=R_{\theta}\cdot z+p,
\label{latticesymmetry}
\end{equation}
where
\[
R_{\theta}\equiv \left(\begin{array}{lr}\cos\,\theta&-\sin\,\theta\\\sin\,\theta&\cos\,\theta\end{array}\right),\,\,\,\,p=\left(\begin{array}{c}p_x\\p_y\end{array}\right).
\]

It thus follows that the group composition is the semi-direct product
\[
(\theta_2,p_2)\cdot(\theta_1,p_1)=(\theta_1+\theta_2,R_{\theta_2}p_1+p_2).
\]

We will denote by $\Sigma$ the following subgroup of $\SEtwo$, which represents the symmetry of a regular square lattice
\[
\Sigma=\left\{\,(\theta,p)\in\SEtwo\,\,|\,\,\theta=\frac{\ell\pi}{2} \mbox{\rm (mod}\,\,2\pi\mbox{\rm )},\,p=\left(\begin{array}{c}m\\n\end{array}\right),\,\,\ell,m,n\in\mathbb{Z}\,\right\}
\]

Let $X$ be a Banach space, and suppose 
\begin{equation}
a:\SEtwo\rightarrow GL(X)
\label{aaction}
\end{equation}
is a faithful and isometric representation of $\SEtwo$ in the space of bounded, invertible linear operators on $X$.  For example, if $X$ is a space of functions with planar domain, a typical $\SEtwo$ action is given by
\[
(a(\gamma) u)(z)=u(\gamma^{-1}\cdot z),\;\;\;\;\gamma=(\theta,p)\in\SEtwo.
\]

We consider an autonomous differential equation on $X$ of the form
\begin{equation}
u_t=f(u)+\varepsilon\,g(u,\varepsilon),
\label{basicpde}
\end{equation}
where $\varepsilon\geq 0$ is a small parameter, $f$ and $g$ satisfy the usual conditions \cite{Henry} guaranteeing that $(\ref{basicpde})$ generates a smooth local semiflow $\Phi_{t,\varepsilon}$ on $X$, and $g$ is bounded.

We will assume certain compatibility conditions between the evolution equation (\ref{basicpde}) and the actions of $\SEtwo$ and $\Sigma$ on $X$ via the representation (\ref{aaction}).  Specifically, we assume
\begin{hyp}
\begin{equation}
\Phi_{t,0}(a(\gamma)u)=a(\gamma)\Phi_{t,0}(u),\;\;\;\forall\, 
u\in X, \gamma\in\SEtwo,
t > 0,
\end{equation}
and for $\varepsilon > 0$, we have
\begin{equation}
\Phi_{t,\varepsilon}(a(\gamma)u)=a(\gamma)\Phi_{t,\varepsilon}(u),
\;\;\;\forall\, u\in X,
\forall\,t > 0 \Longleftrightarrow \gamma\in\Sigma.
\end{equation}
\label{h1}
\end{hyp}
Hypothesis \ref{h1} states that $f$ in (\ref{basicpde}) is 
$\SEtwo-$equivariant,
but $g$ is only $\Sigma-$equivariant.

We want to investigate the effects of the symmetry-breaking perturbation on a normally hyperbolic relative equilibrium to (\ref{basicpde}).  Therefore, we will
postulate that such a solution exists:
\begin{hyp}{\em (Existence of normally hyperbolic relative equilibrium)}
There exists $u^*\in X$ and $\xi^*$ in the Lie algebra of $\SEtwo$
such that
\[
\Phi_{t,0}(u^*)=a(\gamma(t))u^*
\]
holds for all $t$,
where $\gamma(t)=e^{\xi^* t}$.
We also assume that 
the set $\{\lambda\in \mathbb{C}\,|\,|\lambda |
\geq 1\}$ is a spectral set for the linearization
$a(e^{-\xi^*})D\Phi_{1,0}(u^*)$ with projection $P_*$
such that the generalized eigenspace $\mbox{\rm range}(P_*)$ is
three-dimensional (corresponding to the symmetry eigenvalues).
\label{h2}
\end{hyp}

In order to simplify the discussion, we will assume that the isotropy subgroup of $u^*$ in Hypothesis \ref{h2} is trivial, i.e.
\[
(\theta,p)\cdot u^*=u^*\Longleftrightarrow (\theta,p)=(0,0).
\]

Assuming that all other hypothesis of the center manifold theorem of \cite{sandetal2in96,sandetalin98} are satisfied, then for $\varepsilon$ small enough, the dynamics of (\ref{basicpde}) near the relative equilibrium reduce to the following ODE system on the center bundle $Y\cong\Tthree$, (where $\Tthree$ is the three-torus)

\begin{equation}
\begin{array}{rcl}
\dot{\Psi}&=&R_{\varphi}(V+\varepsilon\,F^\Psi(\Psi,\varphi,\varepsilon))\\&&\\
\dot{\varphi}&=&\omega+\varepsilon\,F^{\varphi}(\Psi,\varphi,\varepsilon)
\end{array}
\label{cmeqs1}
\end{equation}

where ${\displaystyle \left(\begin{array}{c}\Psi\\\varphi\end{array}\right)=\left(\begin{array}{c}\psi_1\\\psi_2\\\varphi\end{array}\right)\in\Tthree}$, $V\in\mathbb{R}^2$ is a constant vector, $\omega\geq 0$ is a constant real number, and $\varepsilon\geq 0$ is a small parameter.  The functions $F^\Psi$ and $F^{\varphi}$ are assumed to be smooth enough for our purposes, are $2\pi$-periodic in $\psi_1$, $\psi_2$ and in $\varphi$, and have the following symmetry property
\begin{equation}
F^{\Psi,\varphi}\left(-J\Psi,\varphi+\frac{\pi}{2},\varepsilon\right)=
F^{\Psi,\varphi}\left(\Psi,\varphi,\varepsilon\right),
\,\,\,\,\,\,\,\,\,\,\,\forall \left(\begin{array}{c}\Psi\\\varphi\end{array}\right)\in
\Tthree,\, 0\leq\varepsilon\ll 1,
\label{symprop}
\end{equation}
where
\begin{equation}
J=\left(\begin{array}{rc}0&1\\-1&0\end{array}\right)=R_{-\frac{\pi}{2}}.
\label{Jdef}
\end{equation}
As a consequence, if $(\Psi(t),\varphi(t))$ is a solution of (\ref{cmeqs1}) then so is ${\displaystyle\left(J\Psi(t),\varphi(t)-\frac{\pi}{2}\right)}$.



\Section{Averaging and invariant manifolds}

The main theoretical tool underlying the analysis in the later sections is an approach which is presented in \cite{Ha80} to prove the existence (and asymptotic stability) of invariant manifolds in systems of ordinary differential equations of the form
\arraystart
\begin{equation}
\begin{array}{lll}
\dot{x}&=&\varepsilon A x + \varepsilon F(\phi,x,\varepsilon)\\
\dot{\phi}&=&\omega+\varepsilon\Theta(\phi,x,\varepsilon),
\end{array}
\label{Ha80:eq24}
\end{equation}
\arrayfinish
where $x\in {\bf R}^n$, $\phi\in {\bf R}^m$,
$\omega\in {\bf R}^m$ is a constant vector, 
$A$ is an $n$ by $n$ matrix
of constants of block diagonal form
\[
A=\left(\begin{array}{c|c}
A_s&0\\
\hline
0&A_u\end{array}\right),
\]
where $A_s$ has all eigenvalues with negative real part and $A_u$
has all eigenvalues with positive real part.

The main result we will need is Theorem 2.3, \S VII.2 of \cite{Ha80} which we state here for convenience:
\begin{thm} (\cite[\S VII.2, Theorem 2.3]{Ha80})\newline
Consider the system (\ref{Ha80:eq24}), and
define
\[
\Omega(r_0,\varepsilon_0)=\{(x,\varepsilon):|x| < r_0, 0\leq\varepsilon
\leq\varepsilon_0\}.
\]
Suppose that $\Theta$,
$F$ satisfy the following conditions:
\begin{enumerate}
\item[(C1)] 
$\Theta$ and $F$ are continuous and bounded in ${\bf R}^m\times
\Omega(r_0,\varepsilon_0)$.
\item[(C2)] $\Theta$ and $F$ are Lipschitz in $\phi$
in ${\bf R}^m\times
\Omega(r,\varepsilon)$,
with Lipschitz
constants 
$L^{\Theta}_{\phi}(r,\varepsilon)$ and
$L^{F}_{\phi}(r,\varepsilon)$
respectively, where  
$L^{\Theta}_{\phi}(r,\varepsilon)$, $L^{F}_{\phi}(r,\varepsilon)$
are
continuous and non-decreasing for $0\leq r\leq r_0$, $0\leq\varepsilon
\leq\varepsilon_0$, and $L^{F}_{\phi}(0,0)=0$.
\item[(C3)] $\Theta$ and $F$ are Lipschitz in $x$ with Lipschitz
constants 
$L^{\Theta}_{x}(r,\varepsilon)$ and $L^{F}_{x}(r,\varepsilon)$
respectively, in ${\bf R}^m\times
\Omega(r,\varepsilon)$, where 
$L^{\Theta}_{x}(r,\varepsilon)$, $L^{F}_{x}(r,\varepsilon)$
are
continuous and non-decreasing for $0\leq r\leq r_0$, $0\leq\varepsilon
\leq\varepsilon_0$, and $L^{F}_{x}(0,0)=0$.
\item[(C4)] The function $|F(\phi,0,\varepsilon)|$ is bounded by $N(\varepsilon)$
for $\phi\in {\bf R}^m$, $0\leq\varepsilon\leq\varepsilon_0$, where
$N(\varepsilon)$ is continuous and nondecreasing for
$0\leq\varepsilon\leq\varepsilon_0$, and $N(0)=0$.
\end{enumerate}
Furthermore, 
since the matrix $A$ has all eigenvalues bounded away from the imaginary
axis, there exists 
$\alpha>0$, $K>0$
such that 
for any real number $\tau$, 
\[
\begin{array}{lll}
|e^{(t-\tau)\varepsilon A_s} |&\leq&
Ke^{-\varepsilon\alpha (t-\tau)},\;\;\;\;\;t\geq\tau,
\\
|e^{(t-\tau)\varepsilon A_u}|&\leq&
Ke^{-\varepsilon\alpha (t-\tau)},\;\;\;\;\;t\leq\tau.
\end{array}
\]
Suppose that $\alpha-\limsup_{\varepsilon\rightarrow 0}
L^{\Theta}_{\phi}(0,\varepsilon) > 0$.  

Then the following conclusions are true for the system (\ref{Ha80:eq24}):
there exist $\varepsilon_1 > 0$ and continuous functions
$D(\varepsilon)$, $\Delta(\varepsilon)$, $0 < \varepsilon \leq \varepsilon_1$,
which approach zero as $\varepsilon\rightarrow 0$, and a continuous function
\[
\sigma : {\bf R}^m\times (0,\varepsilon_1]\longmapsto {\bf R}^n
\]
with $\sigma(\phi,\varepsilon)$
bounded by $D(\varepsilon)$, Lipschitz in $\phi$ with Lipschitz constant
$\Delta(\varepsilon)$, such that the set
\[
S_{\varepsilon} = \{(\phi,x) : x = \sigma(\phi,\varepsilon), \phi\in {\bf R}^m\}
\]
is an 
invariant set for the system (\ref{Ha80:eq24}).
If the functions $\Theta$ and $F$ in (\ref{Ha80:eq24}) are multiply periodic
in $\phi$ with period vector $(T_1,\ldots,T_m)$, then $\sigma(\phi,\varepsilon)$ 
is also multiply
periodic in $\phi$ with period vector $(T_1,\ldots,T_m)$.
Finally, the asymptotic stability of the invariant set $S_{\varepsilon}$ is
the same as that of the trivial equilibrium point of the equation
$\dot{y}=Ay$.
\label{Ha80:thm23}
\end{thm}

In the next sections, we will cast equations (\ref{cmeqs1}) in forms suitable for application of Theorem \ref{Ha80:thm23}, and prove the existence of hyperbolic invariant manifolds in these equations.  We will then give physical interpretations for these invariant manifolds.

\Section{Perturbed rotating waves}

When $\omega = \omega_0 > 0$ in (\ref{cmeqs1}), then the right-hand side of the $\dot{\varphi}$ equation is strictly positive for all $\varepsilon > 0$ small enough.  Consequently we may rescale time along orbits of (\ref{cmeqs1}) so that the resulting equations are of the form
\begin{equation}
\begin{array}{rcl}
\dot{\Psi}&=&R_{\varphi}(W+\varepsilon H(\Psi,\varphi,\varepsilon))\\&&\\
\dot{\varphi}&=&1
\end{array}
\label{pertrot}
\end{equation}
where $W$ is a constant vector and $H:\Tthree\times\mathbb{R}\longrightarrow\mathbb{R}^2$ is smooth and satisfies the symmetry properties (\ref{symprop}).
A simple change of variables transforms (\ref{pertrot}) into
\begin{equation}
\begin{array}{rcl}
\dot{\Psi}&=&\varepsilon\,R_{\varphi}G(\Psi,\varphi,\varepsilon)\\&&\\
\dot{\varphi}&=&1
\end{array}
\label{pertrot2}
\end{equation}
where $G(\Psi,\varphi,\varepsilon)=H(\Psi+JR_{\varphi}W,\varphi,\varepsilon)$ satisfies the symmetry properties (\ref{symprop}).
When $\varepsilon=0$ in (\ref{pertrot2}), all solution curves $(\Psi(t),\varphi(t))=(\Psi_0,\varphi_0+t)$ are $2\pi$-periodic in time on the three-torus.  
In physical space, these motions represent a spiral wave which is rigidly and uniformly rotating around a fixed point in the plane.  The initial conditions $\Psi_0$ and $\varphi_0$ for any solution of (\ref{pertrot2}) can be thought of as parametrizing all possible centers of rotation in the plane for the spiral wave (modulo the unit square), and initial phase of rotation of the spiral wave, respectively.

Our main result in this section will be to show that, under certain hypotheses on $G$, when $\varepsilon >0$, (\ref{pertrot2}) can admit hyperbolic periodic solutions, and hyperbolic invariant two-tori, where one of the frequencies of the two-tori is $O(\varepsilon)$.  

We note using standard results from the theory of averaging that there exists a periodic (in $\varphi$) near-identity change of variables which transforms (\ref{pertrot2}) into
\begin{equation}
\begin{array}{rcl}
\dot{\Psi}&=&\varepsilon {\cal G}(\Psi)+\varepsilon^2 R_{\varphi} U(\Psi,\varphi,\varepsilon)\\&&\\
\dot{\varphi}&=&1
\end{array}
\label{avpertrot}
\end{equation}
where $U(\Psi,\varphi,\varepsilon)$ is smooth, multiply periodic, satisfies (\ref{symprop}), and
\begin{equation}
{\cal G}(\Psi)=\frac{1}{2\pi}\,\int_0^{2\pi}\,R_{\varphi} G\left(\Psi,\varphi,0\right)\,d\,{\varphi}.
\label{CalGdef}
\end{equation}
A simple computation shows that 
\[
{\cal G}(J\Psi)=J{\cal G}(\Psi).
\]
The ${\mathbb Z}_4$-equivariant planar ordinary differential equation 
\begin{equation}
\dot{\Psi}={\cal G}(\Psi)
\label{z4eq}
\end{equation}
has an equilibrium at $\Psi=0$, and any non-trivial equilibria occur in conjugate sets $\{J^k\Psi^*,k=0,1,2,3\}$.  Also, if $\{\Psi(t)\,|\,0\leq t\leq T\}$ is a $T$-periodic orbit for (\ref{z4eq}), then either $\{J^k\Psi(t)\}$ is a distinct periodic orbit for $k=1,2,3$, or $\{J\Psi(t)\}$ coincides with the orbit $\{\Psi(t)\}$, in which case the periodic solution has one of the spatial-temporal symmetries
\begin{equation}
\Psi(t-T/4)=\pm J\Psi(t).
\label{spatiotemporalsym}
\end{equation}

We are now ready to state the main result of this section, whose proof is given in Appendix A.
\begin{thm}
Consider the ${\mathbb Z}_4$-equivariant ODE (\ref{z4eq}).
\begin{enumerate}
\item[(i)] Suppose (\ref{z4eq}) has an equilibrium point $\Psi^*$ which is linearly stable (resp. unstable).  Then for all $\varepsilon >0$ small enough, (\ref{avpertrot}) has linearly stable (resp. unstable) $2\pi$-periodic orbits of the form
\[
\Psi=J^k(\Psi^*+ \sqrt{\varepsilon}\sigma_{\Psi^*}(\varphi+k\pi/2,\varepsilon)),\,\, k=0,1,2,3, 
\]
where $\sigma_{\Psi^*}$ approaches $0$ as $\varepsilon\rightarrow 0$.  Furthermore,
if $\Psi^*=0$, then the periodic orbit is such that
\begin{equation}
\sigma_0(\varphi-\pi/2,\varepsilon)=J\sigma_0(\varphi,\varepsilon).
\label{stsp}
\end{equation}
\item[(ii)] Suppose $\{\Psi^*(t)\,|\,0\leq t\leq T\}$ is a non-trivial $T$-periodic orbit for (\ref{z4eq}), and that this orbit is linearly stable (resp. unstable).  Then for all $\varepsilon > 0$ small enough, (\ref{avpertrot}) has linearly stable (resp. unstable) invariant two-tori of the form
\[
\begin{array}{ll}
\Psi=J^k(\Psi^*(\theta+kT/4)+\sqrt{\varepsilon}\,\Sigma_{\Psi^*}(\theta+kT/4,\varphi+k\pi/2,\varepsilon)),&k=0,1,2,3,\\
&(\theta,\varphi)\in [0,T]\times [0,2\pi], 
\end{array}
\]
where $\Sigma_{\Psi^*}$ approaches $0$ as $\varepsilon\rightarrow 0$.  Furthermore, If $\Psi^*$ satisfies the spatio-temporal symmetry property (\ref{spatiotemporalsym}) 
then the invariant torus is such that
\begin{equation}
\Sigma_{\Psi^*}(\theta-T/4,\varphi\mp\pi/2,\varepsilon)=\pm J\Sigma_{\Psi^*}(\theta,\varphi,\varepsilon).
\label{stsp2}
\end{equation}

\end{enumerate}
\label{thmrotwave}
\end{thm}

\noindent
{\bf Interpretation} 

\vspace*{0.1in}
\noindent
As mentioned at the beginning of this section, when $\varepsilon=0$ the dynamics of (\ref{pertrot2}) are particularly simple: all solutions 
are $2\pi$-periodic in time on the three-torus.  
In physical space, these motions represent a spiral wave which is rigidly and uniformly rotating around a fixed point in the plane.  The initial conditions $\Psi_0$ and $\varphi_0$ for any solution of (\ref{pertrot2}) can be thought of as parametrizing all possible centers of rotation in the plane for the spiral wave (modulo the unit square), and initial phase of rotation of the spiral wave, respectively.  Thus in the fully Euclidean case, there are no preferred centers of rotation.  All points in physical space can be centers of rotation provided the proper initial conditions are chosen.

Theorem \ref{thmrotwave} characterizes the existence of (generically isolated) periodic solutions and invariant two-tori for (\ref{pertrot2}) and the local 
structure of the dynamics around these isolated periodic solutions and tori for $\varepsilon > 0$.  In physical space, these hyperbolic periodic solutions and two-tori correspond to 
``anchored'' rotating waves or meandering waves (in the asymptotically stable case), or ``repelling'' rotating waves or meandering waves (in the unstable case) \cite{BLM,LW}. 
Therefore, whereas rotating spiral waves in the fully Euclidean case have no preferred center of rotation, in the perturbed case, we expect generically that certain points in physical space (called {\em AR-points} herein) will act as local attractors or repellants for the spiral core.   To any such AR-point in physical space, there correspond in fact infinitely many conjugate AR-points related to each other by the lattice symmetry (\ref{latticesymmetry}).  Perhaps counter-intuitively, a given AR-point need not coincide with the physical location of localized heterogeneities (i.e. lattice points), as was shown in \cite{BLM}.  Although if it does, then the solution possesses residual spatio-temporal symmetry.  The phenomenon is illustrated with numerical simulations in Section 5 below.   

To summarize, Theorem \ref{thmrotwave} states that the effects of lattice symmetry-breaking on a rotating spiral wave in a Euclidean system will be to create hyperbolic rotating waves and/or hyperbolic meandering waves.  In the latter case, one of the meander frequencies is very small (of the order of the ``size'' of the perturbation).   In certain cases (where the AR-point of the wave coincides with a point of the lattice), the rotating wave or meandering wave has spatial temporal symmetry: for the rotating wave, rotation in physical space by $\pi/2$ around the lattice point is the same as time-evolution by a quarter period.  In the meandering case, the full meander pattern will have symmetry by rotation of $\pi/2$ around the lattice point.  All of these possibilities are illustrated in the numerical examples in Section 6.

\Section{Perturbed travelling waves}

When $\varepsilon=0$ in (\ref{cmeqs1}), these equations represent travelling waves when $\omega=0$.  
In physical space, these travelling waves manifest themselves as a spiral wave that has completely unwound, with retracting tip \cite{JW,MZ}.

On the three torus, the corresponding solution curves are of the form
\[
(\Psi(t),\varphi(t))=(R_{\varphi_0}Vt+\Psi_0,\varphi_0),
\]
where $\Psi_0$ and $\varphi_0$ are constants.  
One can view these solutions as inducing a partition of the three-torus into a one-parameter family (parametrized by $\varphi_0$) of invariant two-tori, and on each of these invariant two-tori, the flow is either dense or periodic, depending on whether or not the components of the vector $R_{\varphi_0}V$ are linearly independent over the rationals.  We now wish to give a characterization of what happens to this picture when $\varepsilon>0$ is small.  We have the following, which is proved in Appendix B:

\begin{thm}
For the function $F^{\varphi}$ and the constant vector $V$ in (\ref{cmeqs1}), define
\begin{equation}
M(\varphi)=\frac{1}{4\pi^2}\iint_{\Ttwo}\,F^{\varphi}(\Psi,\varphi,0)\,d\Psi\,.
\label{Mdef}
\end{equation}
Suppose $\varphi^*$ is such that
\begin{enumerate}
\item[(i)] $M(\varphi^*)=0$,
\item[(ii)] $\mu=M'(\varphi^*)\neq 0$.
\end{enumerate}
Further suppose that the vector ${\displaystyle\left(\begin{array}{c}\alpha\\\beta\end{array}\right)=R_{\varphi^*}V}$ is such that the following Diophantine condition is satisfied for some constants $\sigma>0$ and $K>0$:
\begin{equation}
|\,n_1\alpha+n_2\beta\,|\geq\frac{K}{(|\,n_1\,|+|\,n_2\,|)^{\sigma+2}},\,\,\,\,\forall\,(n_1,n_2)\in\mathbb{Z}^2\setminus\,\{(0,0)\,\}.
\label{diophantine}
\end{equation}
Then for all $\varepsilon>0$ sufficiently small, the system (\ref{cmeqs1}) has a hyperbolic invariant two-torus
\[
\varphi={\cal T}(\Psi,\varepsilon)
\]
such that ${\cal T}(\Psi,\varepsilon)\rightarrow\varphi^*$ as $\varepsilon\rightarrow 0$.  This two-torus is locally asymptotically stable (resp. unstable) if $\mu<0$ (resp. $\mu>0$).
\label{thmtravwave}
\end{thm}

\noindent
{\bf Interpretation}

The existence of stable travelling waves (when the conditions of Theorem \ref{thmtravwave}
are met) is of some interest.  Suppose that $\varphi^*\in\Sone$ is such that
$M(\varphi^*)=0$ and
$\mu=M(\varphi^*)<0$.
Then there exists 
$\varepsilon_0>0$ such that for all $\varepsilon\in (0,\varepsilon_0)$, 
system (\ref{cmeqs1}) has a hyperbolic invariant two-torus
\[
\varphi={\cal T}(\Psi,\varepsilon)
\]
which represents a travelling wave.

\begin{figure}[tbhp]
\centerline{%
\includegraphics[width=3.2in]{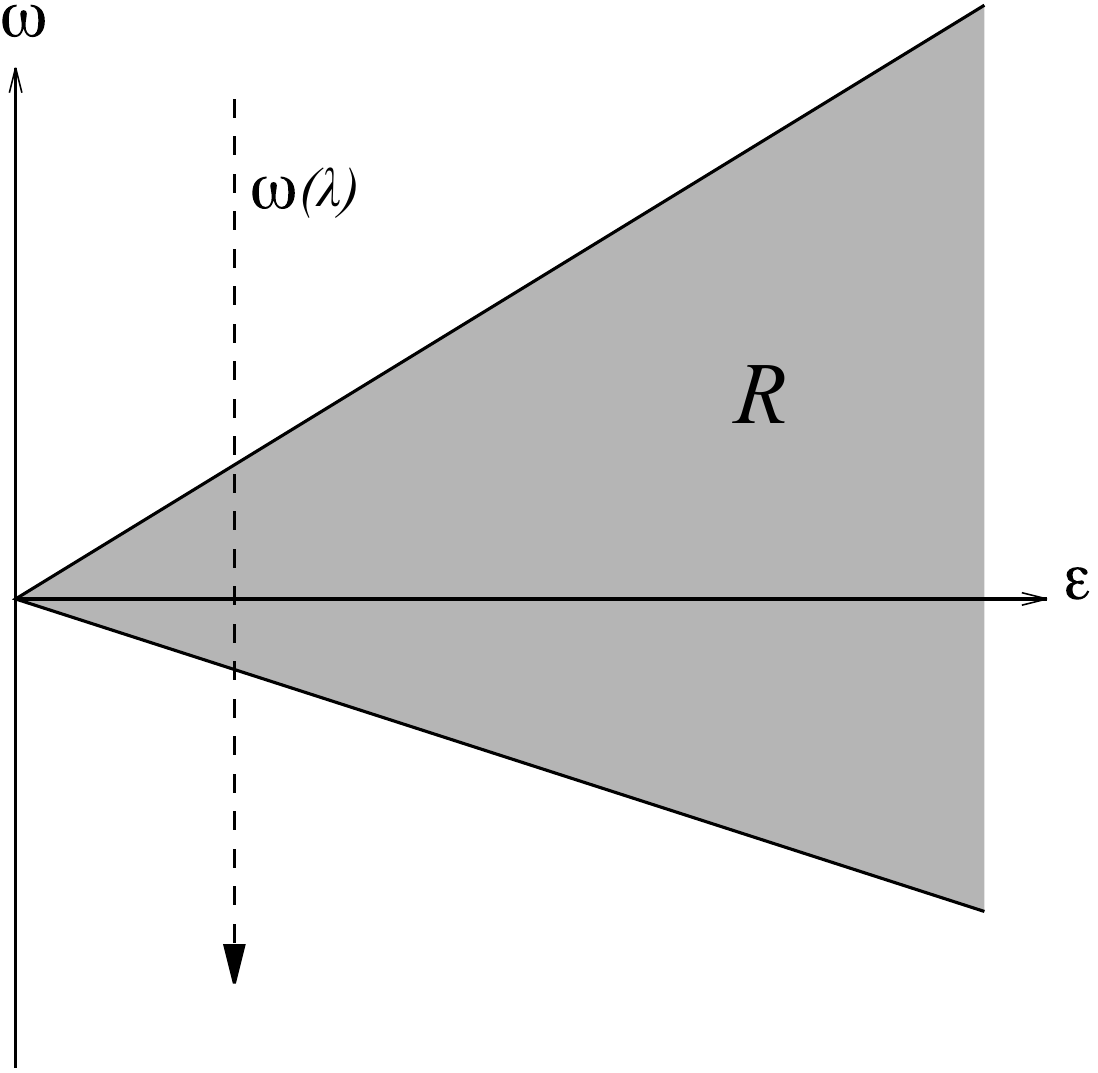}}
\caption{Qualitative sketch of the region $R$ in which
stable travelling waves of system (\ref{cmeqs1}) exist.}
\label{retreg}
\end{figure}

By normal hyperbolicity, we expect that
for all $\varepsilon\in (0,\varepsilon_0)$, there exists
an interval $(-\omega(\varepsilon),\omega(\varepsilon))$ such that 
(\ref{cmeqs1}) has a hyperbolic invariant two-torus
\[
\varphi={\cal T}_{\omega}(\Psi,\varepsilon)
\]
for all $\omega\in (-\omega(\varepsilon),\omega(\varepsilon))$. 
Thus, in the $\varepsilon$ -- $\omega$ parameter space of (\ref{cmeqs1}),
there exists a region $R$ 
which resembles qualitatively Figure \ref{retreg} (see also
Figure 4 of \cite{LambW} for an analogous phenomenon)
with the following properties:
\begin{itemize}
\item outside $R$ where $\omega$ dominates $\varepsilon$, the results of section 4 apply and there are no stable travelling waves for
(\ref{cmeqs1})
\item on the boundary of $R$, stable travelling waves (invariant two-tori) are created via saddle-node bifurcation,
\item in the interior of $R$ where $\varepsilon$ dominates $\omega$, 
Theorem \ref{thmtravwave} applies and there may exist transverse zeros of the function $M(\varphi)$ which correspond to invariant two-tori
$\varphi={\cal T}_{\omega}(\Psi,\varepsilon)$
for (\ref{cmeqs1}). Each such stable two-torus
corresponds to an asymptotically stable travelling wave solution of (\ref{cmeqs1}).
\end{itemize}

Consider now a physical experiment on spiral waves
where Euclidean symmetry is broken by a weak perturbation which preserves lattice symmetry
(i.e $\varepsilon$ is fixed, small and
positive).  Suppose an experimental parameter $\lambda$ is
varied continuously in such a way 
that 
in the $\dot{\varphi}$ equation of (\ref{cmeqs1}), $\omega=\omega(\lambda)$
decreases as illustrated in the dashed line of Figure \ref{retreg}.
The following may be observed in physical space.
As $\omega$ decreases while remaining above the top
boundary of $R$, the
spiral tip will undergo a periodic motion.  However, when $\omega$ crosses the
boundary of $R$, the wave will stop rotating and
start drifting linearly, due to a
saddle-node bifurcation of zeros in the function $M(\varphi)$.  Moreover,
this linear motion will
persist as the experimental
parameter $\lambda$ is further varied while $\omega$ stays
in the interior of $R$.  See also \cite{LeBlanc} for a similar scenario in the case of rotational symmetry-breaking for spiral waves.

The above scenario is reminiscent of the experimentally
observed phenomenon of retracting spiral
waves.  In certain numerical experiments on spiral waves \cite{JW, MZ}, it has been
observed that as an experimental parameter is
varied, the rotation frequency of a
rotating spiral wave decreases, the spiral core unwinds and the
radius of the circular motion grows unbounded, eventually reaching a
state where the spiral is completely unwound, and the tip drifts
linearly (retracts) when the rotation frequency approaches zero.
As the experimental parameter is further varied,
the spiral remains
unwound and the
linearly drifting state persists in an asymptotically stable
manner.

In \cite{AMN99}, it was shown that the Euclidean equivariant equations
(\ref{cmeqs1}) with $\varepsilon=0$ are in partial agreement
with the retracting wave scenario.  
Indeed, suppose that the dashed line in Figure 
\ref{retreg} were drawn a bit more to the left, coinciding with the
$\omega$-axis (i.e. $\varepsilon=0$ for the Euclidean case).  Then, as is shown in
\cite{AMN99}, as $\omega$ decreases, the circular
spiral tip motion grows unbounded like $1/\omega$.  However,
since the intersection of
the region $R$ with the $\omega$-axis consists solely of
the point $\omega=0$,
the analysis predicts that the linearly translating wave exists only
for $\omega=0$ (i.e. it is not asymptotically stable, and does not persist
under further variation of the parameter).  This prediction is
inconsistent with the above-described experiments on retracting
spirals \cite{JW, MZ}.  However, if one considers a situation
where Euclidean symmetry is broken (i.e. $\varepsilon > 0$), then, as 
our analysis reveals,
asymptotically stable travelling waves near $\omega=0$ are
possible.  We note that the PDE model which is integrated in 
\cite{JW, MZ} is fully Euclidean symmetric; however, 
the spatial discretization schemes which are used to perform the
numerical integrations break the symmetry but preserve a lattice symmetry.
We qualify these remarks with two important caveats.  First,
if indeed Euclidean symmetry-breaking (e.g. of a numerical type) were responsible for
the experimentally observed retracting tip phenomenon (we are not claiming that it
is), then the rotation of the
spiral just before the transition to linear translation would be
markedly non-uniform.  As far as we are aware, this non-uniformity is not observed.
Second, as was the case in \cite{AMN99}, the center bundle equations
(\ref{cmeqs1}) predict that as $\omega$ is further varied (beyond
the region $R$), the spiral should start rotating in the
opposite direction.  
This does not appear to be what happens in the experiments:
as the experimental parameter is
further varied, the excitability of the medium is eventually diminished
to an extent where
wave propagation is no longer possible \cite{JW, MZ}.

\Section{Simulations}

\subsection{Anchoring in the perturbed FitzHugh-Nagumo system}

We have performed numerical simulations to illustrate the theoretical results presented in the previous sections.  
We have used a finite-difference scheme ($200\times 200$ spatial grid) with explicit time-stepping to integrate the FitzHugh-Nagumo reaction diffusion PDE
\begin{figure}[b]
\begin{center}
\includegraphics[width=5.3in]{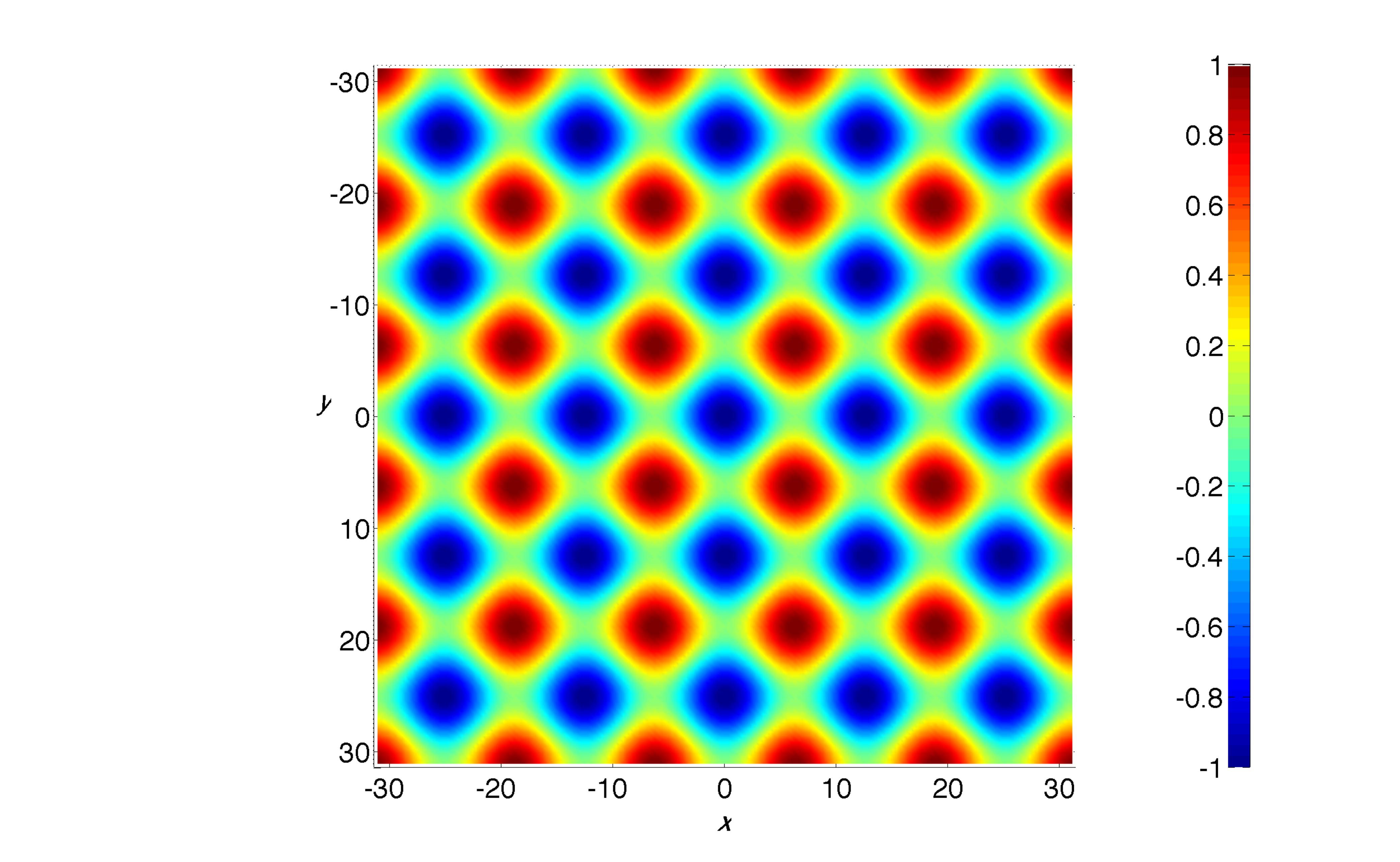}
\caption{First (normalized) spatial harmonic $0.5 (\cos(0.5x)+\cos(0.5y))$ for the inhomogeneous terms $g_1$ and $g_2$ in (\ref{RDPDEsim}).}
\label{inhfield}
\end{center}
\end{figure}
\begin{figure}[htpb]
\begin{center}
\includegraphics[width=3.2in]{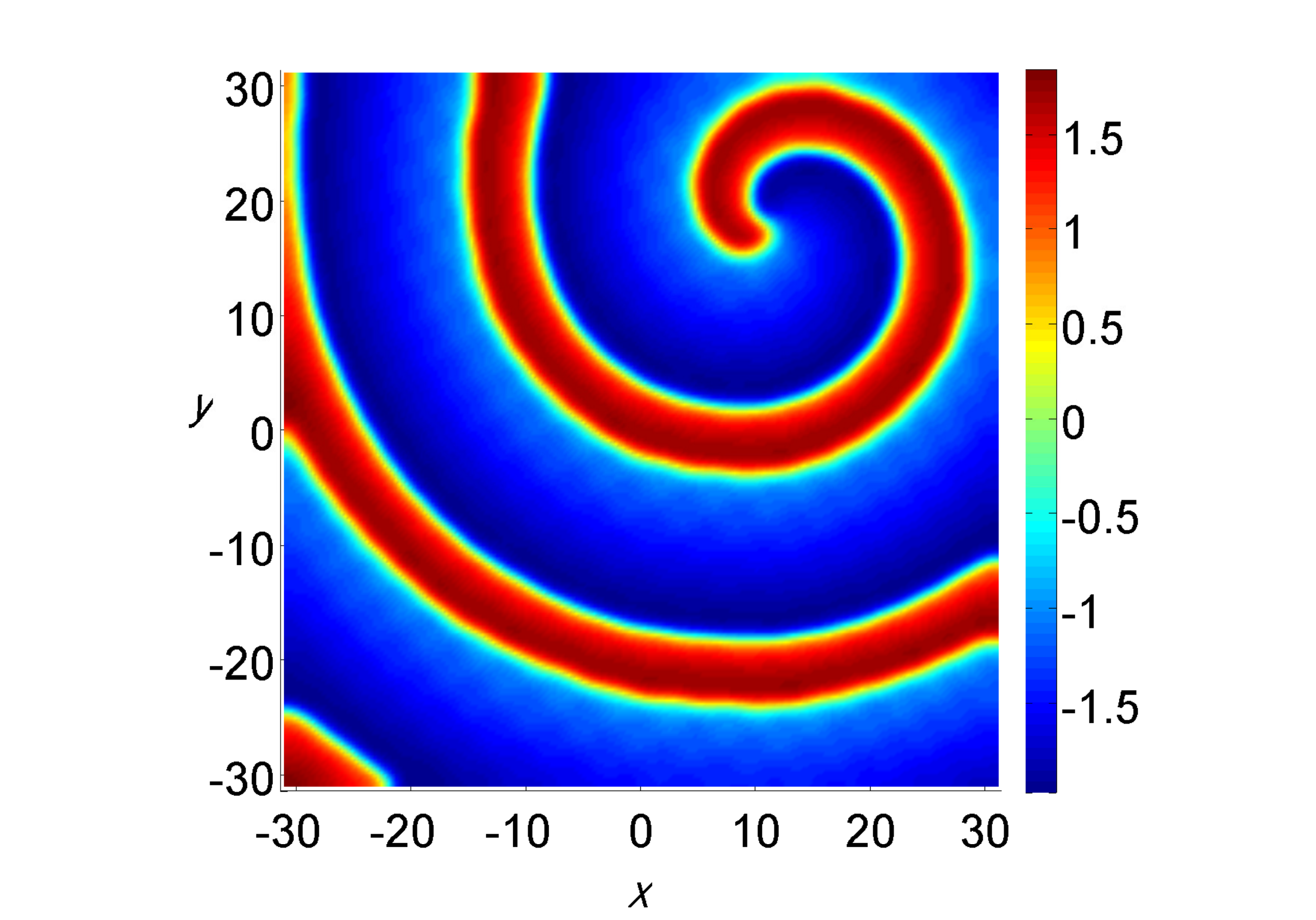}
\includegraphics[width=3.2in]{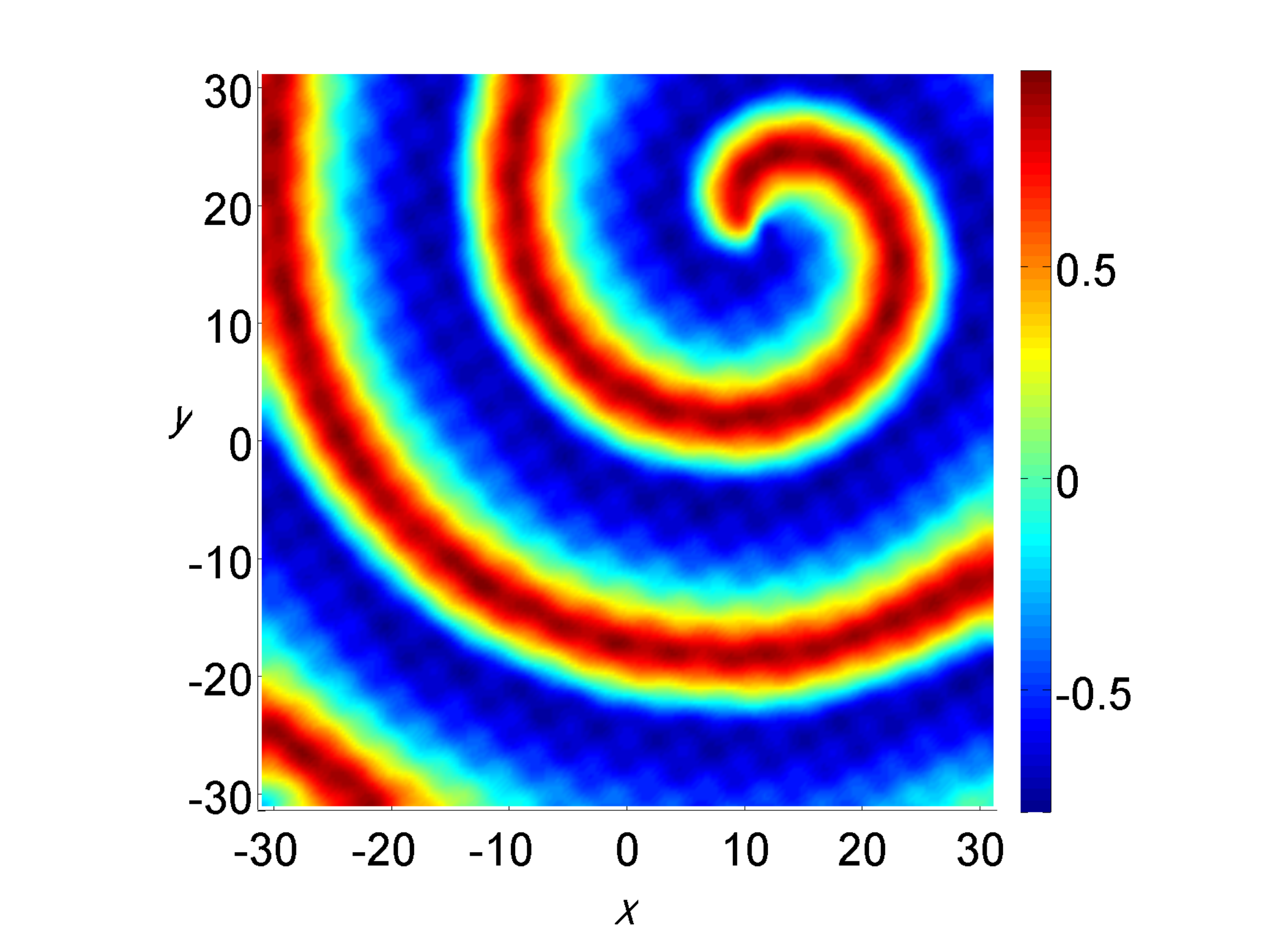}
\caption{Fields $u(x,y)$ (left) and $v(x,y)$ (right) at some instant in time for the integration of (\ref{RDPDEsim}) in the inhomogeneous case.  One can visually observe the background effects of the lattice perturbation on these fields.}
\label{uvfields}
\end{center}
\end{figure}
\begin{figure}[htpb]
\begin{center}
\includegraphics[width=6.5in]{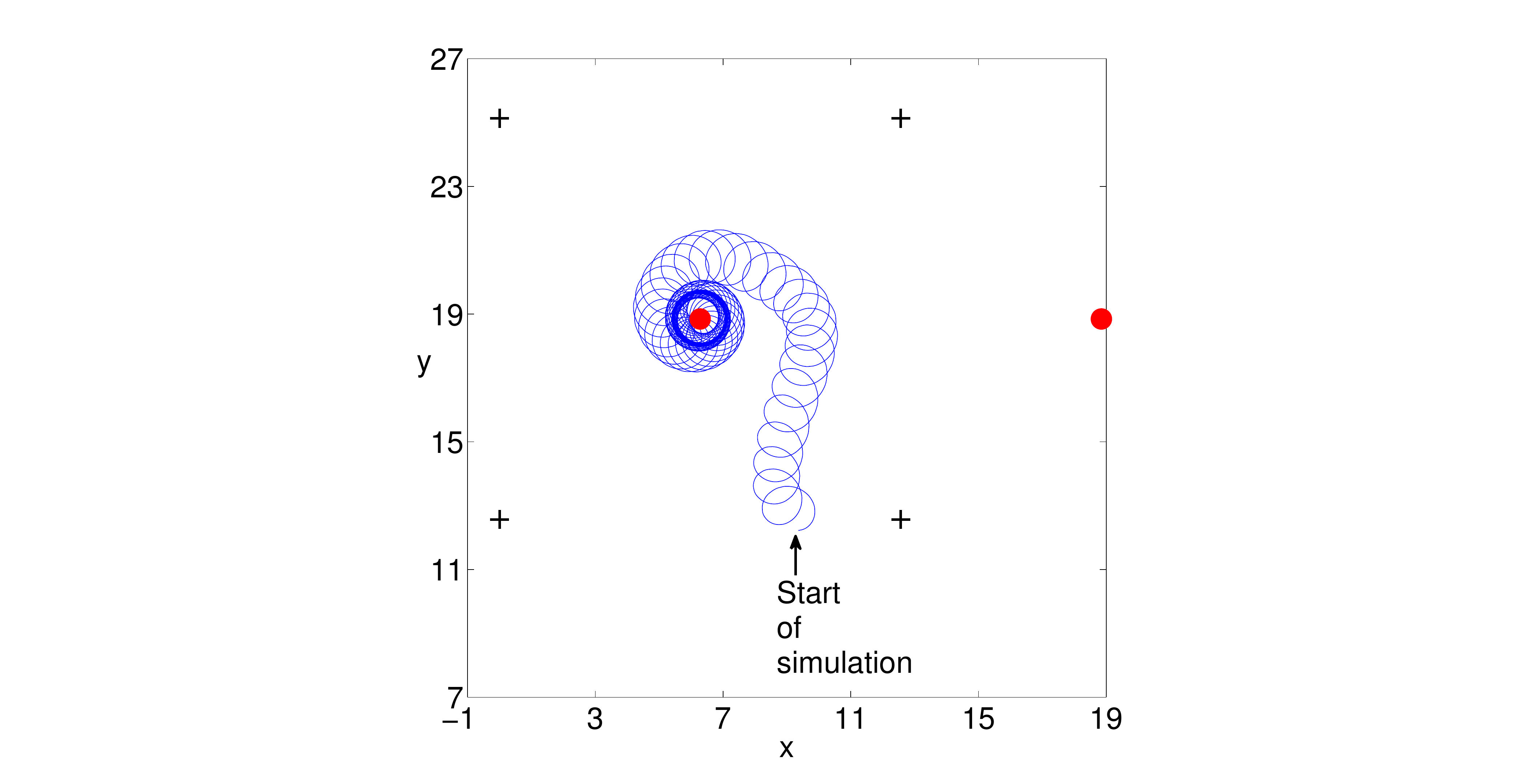}
\end{center}
\caption{Path of the spiral tip for an integration of (\ref{RDPDEsim}) with coefficients (\ref{firstexp}).  The black crosses and the red circles denote the positions of the minima and maxima (respectively) of the first spatial harmonic of the inhomogeneous terms in (\ref{RDPDEsim}).}
\label{anctrans}
\end{figure}
\begin{figure}[htpb]
\begin{center}
\includegraphics[width=6.5in]{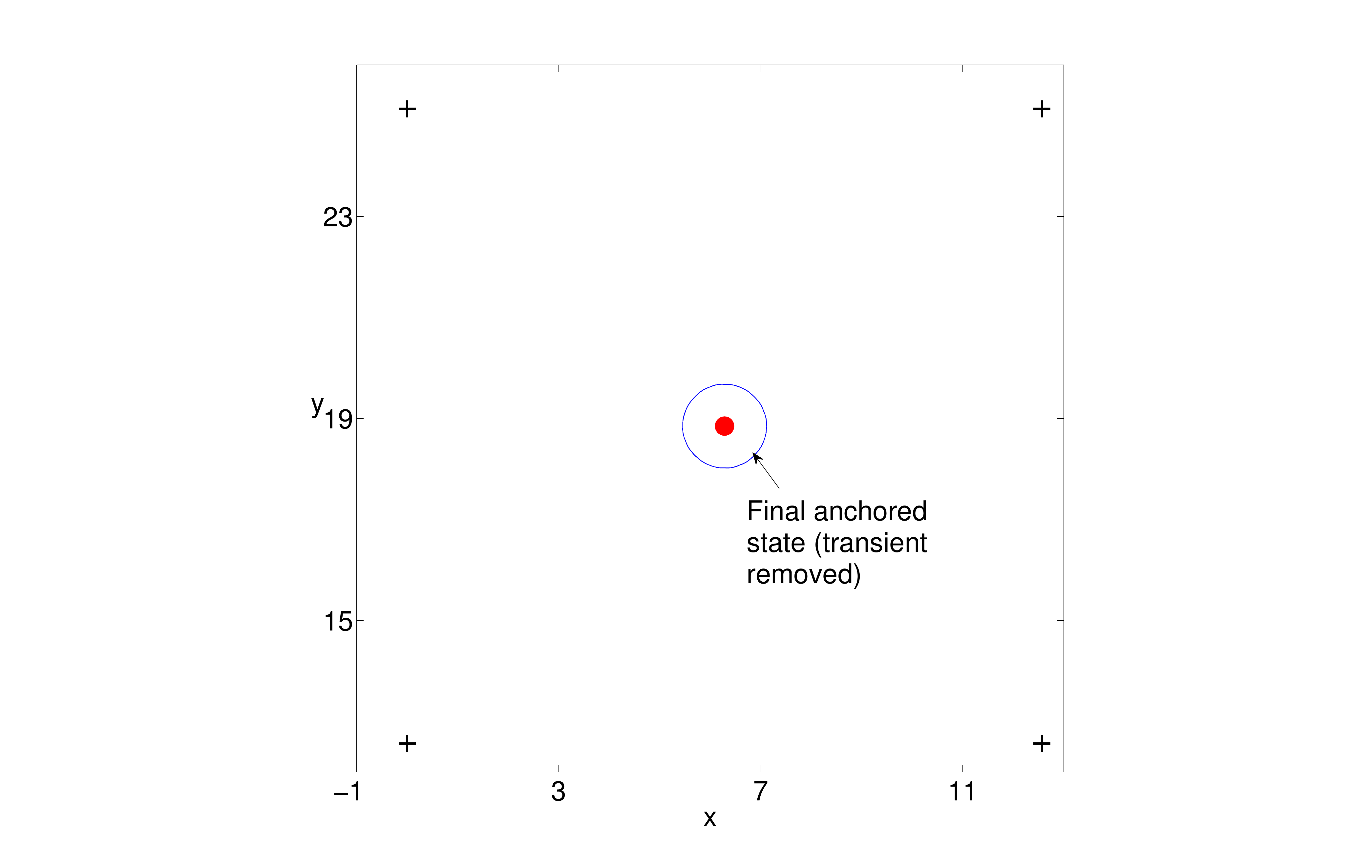}
\caption{Path of the spiral tip for the same integration of (\ref{RDPDEsim}) as in previous figure, with transient removed.  We can see that the final (steady) state is anchored 
around a local maximum of the inhomogeneity field (lattice point).}
\label{ancsteady}
\end{center}
\end{figure}
\begin{figure}[htpb]
\begin{center}
\includegraphics[width=4.5in]{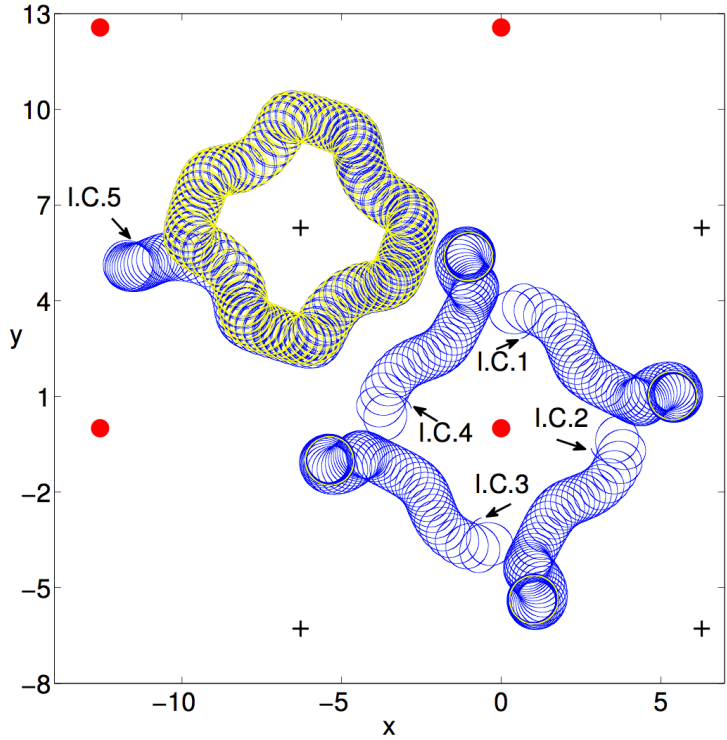}
\caption{Path of the spiral tip for five separate integrations of (\ref{RDPDEsim}) with coefficients (\ref{secondexp}).  
The initial position of the spiral tip for each integration is labeled as I.C.1 through I.C.5.  Initial conditions I.C.1 through I.C.4 are conjugate by rotations of 90 degrees, whereas I.C.5 is not related by symmetry to the other initial conditions.
The final (stable) states are shown in yellow.  The system displays multistability between families of rotating waves which are anchored around non-lattice points, and spatio-symmetric meandering waves, anchored at (dual) lattice points.}
\label{anc4}
\end{center}
\end{figure}
\begin{figure}[htpb]
\begin{center}
\includegraphics[width=4.5in]{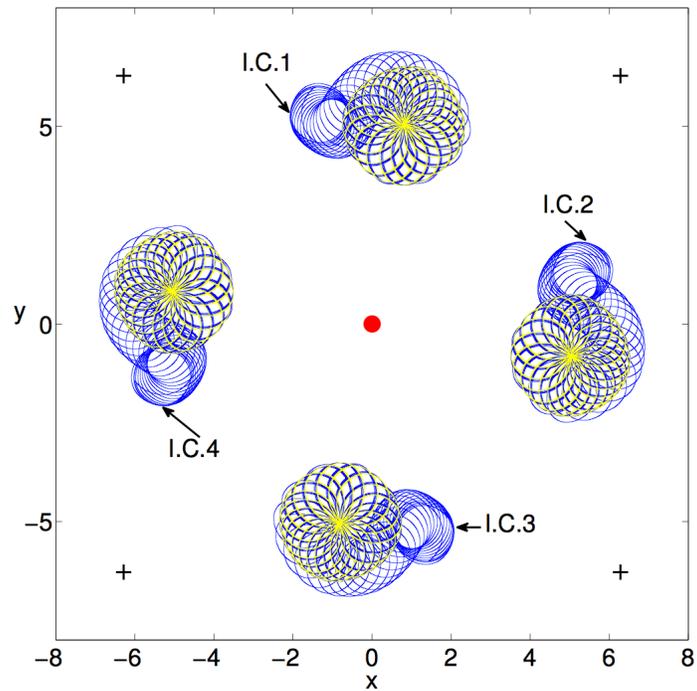}
\caption{Path of the spiral tip for four separate integrations of (\ref{RDPDEsim}) with coefficients (\ref{thirdexp}). The initial position of the spiral tip for each integration is labeled as I.C.1 through I.C.4. which are conjugate by rotations of 90 degrees.  The final (stable) states are shown in yellow, and correspond to meandering waves which are anchored at non-lattice points.}
\label{ff}
\end{center}
\end{figure}
\begin{figure}[htpb]
\begin{center}
\includegraphics[width=6.5in]{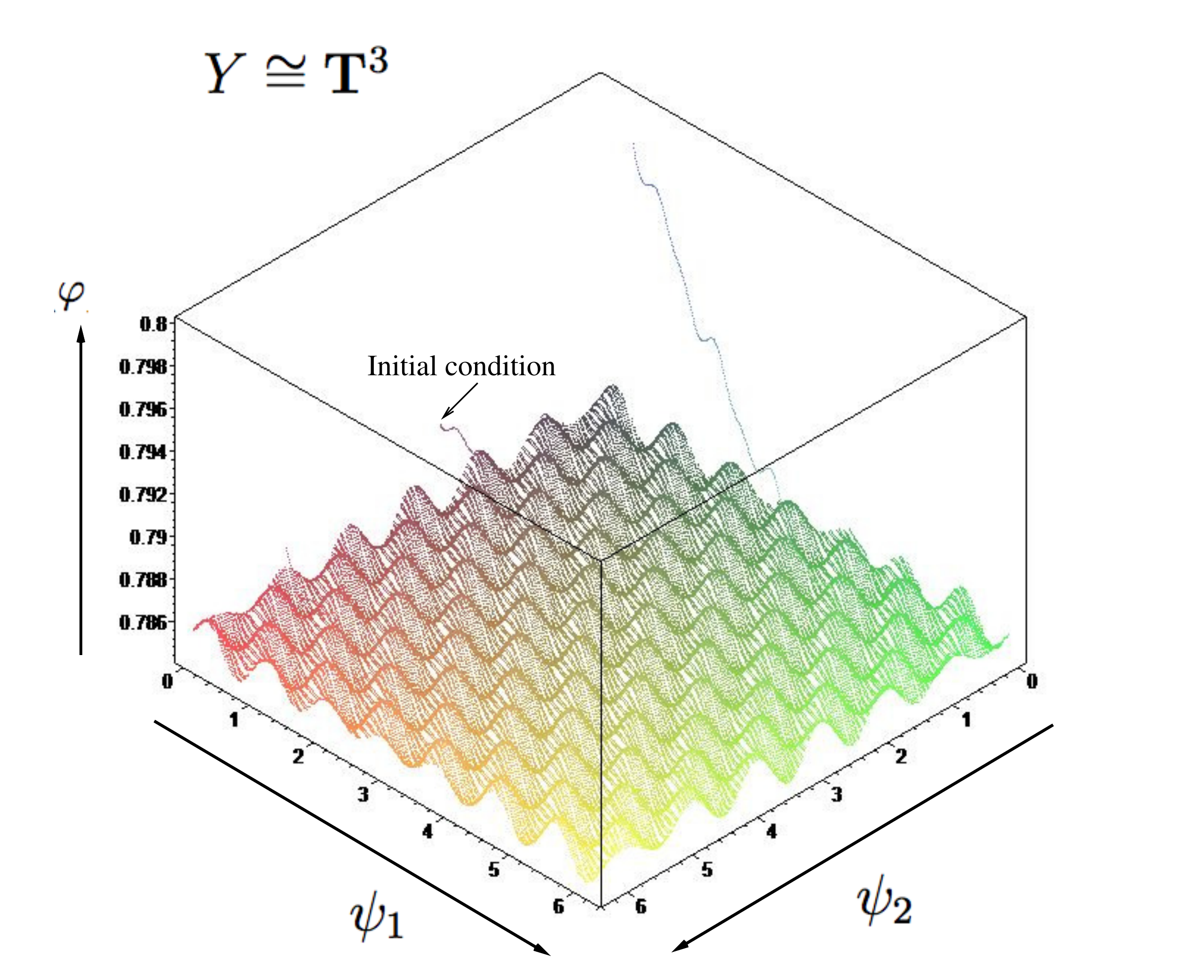}
\caption{Solution curve for an integration of (\ref{cmeqs1}) with data (\ref{simdata}) represented on a portion of the three-torus (opposite faces in $\Psi$ are identified).  After a transient approach, the solution curve fills out densely an invariant surface near the two-torus $\varphi=\pi/4$, as predicted by Theorem \ref{thmtravwave}.}
\label{invs}
\end{center}
\end{figure}

\begin{equation}
\begin{array}{ccl}
{\displaystyle\frac{\partial u}{\partial t}}&=&{\displaystyle\nabla^2\,u+\frac{10}{3}\left(u-\frac{1}{3}u^3-v\right)+\varepsilon\,g_1(x,y)}\\[0.2in]
{\displaystyle\frac{\partial v}{\partial t}}&=&{\displaystyle\frac{3}{10}\left(u+0.6-0.5v\right)+\varepsilon\,g_2(x,y)}
\end{array}
\label{RDPDEsim}
\end{equation}
on the square grid $[-10\pi,10\pi]^2$ with Neumann boundary conditions.  

With the chosen parameters, in the fully Euclidean case $\varepsilon=0$, system (\ref{RDPDEsim}) admits rigidly rotating spiral waves.  We therefore integrate (\ref{RDPDEsim}) with $\varepsilon=0$ long enough so the system eventually achieves a fully-formed, uniformly rigidly rotating spiral wave.  We then use this state as an initial condition for the inhomogeneous case $(\varepsilon\neq 0)$ in further integrations of (\ref{RDPDEsim}).  It should be noted that for all integrations performed in this section, the ``spiral tip'' is defined as the intersection of the $u=0$ and $v=0$ contours (which is computed numerically), but the results would remain qualitatively the same for any other commonly used definition of spiral tip (e.g. maximum of $||\nabla u \times \nabla v||$) \cite{Bar92, Bar94, BKT90, Roth1, Roth2}.

The inhomogeneous terms $g_1$ and $g_2$ in (\ref{RDPDEsim}) are chosen to have the lattice symmetry generated by translations of $4\pi$ along $x$ and $y$ axes, and rotation by 90 degrees around the origin.
Specifically, both $\varepsilon\,g_1$ and $\varepsilon\,g_2$ are expressions of the form
\begin{equation}
\varepsilon\,g_{k}(x,y)=A_{k} + B_k (\cos(0.5x)+\cos(0.5y)) + C_k (\cos(0.5(3x-y))+\cos(0.5(x+3y))),\,\,\,\,\,k=1,2,
\label{lincombcoeffs}
\end{equation}
with the constant coefficients $A_k$, $B_k$ and $C_k$ being changed from one experiment to the next.
The first harmonic $(\cos(0.5x)+\cos(0.5y))$ is plotted in figure \ref{inhfield}.   This setup could represent phenomenologically a spatially periodic (along a lattice) array of reaction inhomogeneities in the medium of propagation of the spiral.

Equations (\ref{RDPDEsim}) are integrated until a steady (anchored) state is observed (either periodic or meandering).  
In figure \ref{uvfields}, we show plots of both the $u$ field and the $v$ field at some intermediate instant in time along the integration.  

Figures \ref{anctrans} and \ref{ancsteady} show the motion in physical space of the spiral tip path (both transients, and steady regimes) for a first choice of coefficients $A_k$, $B_k$ and $C_k$ in (\ref{lincombcoeffs}):
\begin{equation}
\begin{array}{ll}
A_1=0.028&A_2=-0.0044\\
B_1=0.05&B_2=-0.02\\
C_1=0.06&C_2=0.01
\end{array}
\label{firstexp}
\end{equation}
where we see that the spiral gets anchored to the site of a local maximum for the first spatial harmonic of the inhomogeneity field (i.e. anchored to a lattice point).

In figure \ref{anc4}, we show the motion in physical space of the spiral tip path for the following choice of coefficients $A_k$, $B_k$ and $C_k$ in
(\ref{lincombcoeffs})
\begin{equation}
\begin{array}{ll}
A_1=0.016&A_2=0.006\\
B_1=0.05&B_2=-0.0001\\
C_1=0.0001&C_2=0.03
\end{array}
\label{secondexp}
\end{equation}
and for five different integrations using conjugate initial conditions.  
Four of the initial conditions are conjugate by rotations of 90 degrees, whereas the fifth initial condition is non-conjugate.  
As predicted in Section 4, we see that the system displays multistability between conjugate families of stable (anchored) rotating waves, where the point of anchoring is not a lattice point, and meandering waves which are anchored at (dual) lattice points.

Finally, in figure \ref{ff}, we show the motion in physical space of the spiral tip path for the coefficients
\begin{equation}
\begin{array}{ll}
A_1=-0.016&A_2=-0.012\\
B_1=-0.05&B_2=0.0001\\
C_1=-0.0001&C_2=-0.06
\end{array}
\label{thirdexp}
\end{equation}
in (\ref{lincombcoeffs}), with four conjugate initial conditions.  In this case, the stable steady-states are conjugate families of meandering waves anchored at non-lattice points.

\subsection{Invariant surfaces for perturbed travelling waves}

In order to illustrate Theorem \ref{thmtravwave}, we have performed a numerical simulation of the following system of the form (\ref{cmeqs1}) with:
\begin{equation}
\begin{array}{c}
\varepsilon=0.1,\,\,\,\,\,
\omega=0,\,\,\,\,\,
{\displaystyle V=\left(\begin{array}{c}\pi \\ \sqrt{2}\end{array}\right)}\\[0.15in]
F^{\varphi}(\Psi,\varphi,\varepsilon)=2\sin(4\varphi)+\cos(7\psi_1+6\psi_2)+\cos(6\psi_1-7\psi_2)\\[0.15in]
{\displaystyle F^{\Psi}(\Psi,\varphi,\varepsilon)=\left(\begin{array}{c}
\sin(5\varphi)\sin(\psi_1+\psi_2)+\cos(5\varphi)\sin(\psi_1-\psi_2)\\
\cos(2\varphi)(\cos(2\psi_1+3\psi_2)-\cos(3\psi_1-2\psi_2))
\end{array}\right)}.
\end{array}
\label{simdata}
\end{equation}
With this data, the function $M(\varphi)$ in (\ref{Mdef}) becomes
\[
M(\varphi)=2\sin(4\varphi),
\]
and we note in particular that $\varphi^*=\pi/4$ is such that $M(\varphi^*)=0$ and that $\mu=M'(\varphi^*)=-8<0$.  Furthermore, the components of the vector $R_{\varphi^*}V$ are rationally incommensurable, and since the perturbations terms $F^{\varphi}$ and $F^{\Psi}$ are trigonometric polynomials, the diophantine condition (\ref{diophantine}) is satisfied.  Theorem \ref{thmtravwave} then predicts the existence of a locally asymptotically stable invariant surface (two-torus) $\varphi={\cal T}(\Psi,\varepsilon)$ close to $\varphi=\pi/4\approx 0.785$.  Figure \ref{invs} illustrates the result of a numerical simulation of the above system, and we can clearly see the transient approach of the solution to the predicted invariant surface.  It should be noted that this figure represents one integration of the above system for one particular choice of initial condition, and that the trajectory fills out the invariant two-torus densely.

\vspace*{0.25in}
\noindent
{\Large\bf Acknowledgments}

\vspace*{0.2in}
This research is partly supported by the
Natural Sciences and Engineering Research Council of Canada in the
form of a Discovery Grant (VGL), and an Undergraduate Student Research Award (LC).  Funding from the Ontario Graduate Scholarship (OGS) program is also acknowledged (LC).

{\Large\bf\flushleft Appendix}

\appendix

\Section{Proof of Theorem \ref{thmrotwave}}

We first prove item (i) related to equilibria of (\ref{z4eq}).
Suppose $\Psi^*$ is an equilibrium point for (\ref{z4eq}), and write 
\begin{equation}
\Psi=\Psi^*+\sqrt{\varepsilon}\hat{\Psi}.
\label{variation1}
\end{equation}
Then upon dropping the hats, (\ref{avpertrot}) becomes
\begin{equation}
\begin{array}{rcl}
\dot{\Psi}&=&\varepsilon D_{\Psi}{\cal G}(\Psi^*)\Psi+\varepsilon^{\frac{3}{2}} {\cal Q}(\Psi,\varphi,\varepsilon)\\&&\\
\dot{\varphi}&=&1,
\end{array}
\label{avpertrotcentered}
\end{equation}
where ${\cal Q}$ is smooth and periodic in $\varphi$.  
By applying Theorem \ref{Ha80:thm23} in the case $n=2$ and $m=1$, we get that there exists a periodic solution of (\ref{avpertrotcentered}) which is a graph
\[
\Psi=\sigma_{\Psi^*}(\varphi,\varepsilon)
\]
where $\sigma_{\Psi^*}(\varphi,\varepsilon)\rightarrow 0$ as $\varepsilon$ goes to zero, and using (\ref{variation1}) we get a periodic solution of (\ref{avpertrot})
\begin{equation}
\Psi=\Psi^*+\sqrt{\varepsilon}\sigma_{\Psi^*}(\varphi,\varepsilon).
\label{persol1}
\end{equation}
The stability of this periodic solution is the same as that of the equilibrium point $\Psi^*$ in the planar system (\ref{z4eq}).

For the conjugate equilibrium $J\Psi^*$ of (\ref{z4eq}), we replace (\ref{variation1}) with
\begin{equation}
\Psi=J(\Psi^*+\sqrt{\varepsilon}\hat{\Psi})
\label{variation2}
\end{equation}
and set $\phi=\varphi+\pi/2$ in (\ref{avpertrot}).  Because of the equivariance properties (\ref{symprop}), upon dropping the hats, we get
\begin{equation}
\begin{array}{rcl}
\dot{\Psi}&=&\varepsilon D_{\Psi}{\cal G}(\Psi^*)\Psi+\varepsilon^{\frac{3}{2}} {\cal Q}(\Psi,\phi,\varepsilon)\\&&\\
\dot{\phi}&=&1,
\end{array}
\label{avpertrotcentered2}
\end{equation}
which is essentially identical to (\ref{avpertrotcentered}).
So we have a periodic solution
\begin{equation}
\Psi=J\Psi^*+J\sqrt{\varepsilon}\sigma_{\Psi^*}(\phi,\varepsilon)=J\Psi^*+J\sqrt{\varepsilon}\sigma_{\Psi^*}(\varphi+\pi/2,\varepsilon)
\label{persol2}
\end{equation}
for (\ref{avpertrot}).  If $\Psi^*\neq 0$, then (\ref{persol1}) and (\ref{persol2}) are distinct periodic solutions of (\ref{avpertrot}).  However, if $\Psi^*=0$, then 
(\ref{persol1}) and (\ref{persol2}) are the same solution, and we thus get
\[
\sigma_0(\varphi,\varepsilon)=J\,\sigma_0(\varphi+\pi/2,\varepsilon).
\]
which establishes (\ref{stsp}).

Now, we will prove item (ii) related to periodic orbits of (\ref{z4eq}).
Let $\Psi^*(t)$ be a $T$-periodic solution of (\ref{z4eq}).  Then it is well-known \cite{Ha80,Hartman} that near this periodic orbit, there exist local coordinates $(\rho,\theta)$ defined by
\begin{equation}
\Psi=\Psi^*(\theta)+\rho\,J\,\dot{\Psi}^*(\theta)
\label{cov1}
\end{equation}
which transforms (\ref{avpertrot}) into
\begin{equation}
\begin{array}{rcl}
\dot{\rho}&=&\varepsilon\,A(\rho,\theta)\rho+\varepsilon^2 Q(\rho,\theta,\varphi,\varepsilon)\\&&\\
\dot{\theta}&=&\varepsilon\,(1+B(\rho,\theta)\rho)+\varepsilon^2 R(\rho,\theta,\varphi,\varepsilon)\\&&\\
\dot{\varphi}&=&1,
\end{array}
\label{pertavg1}
\end{equation}
where all functions are smooth, $T$-periodic in $\theta$ and $2\pi$-periodic in $\varphi$, and where the function $A(\rho,\theta)$ is such that
\begin{equation}
A(0,\theta)=\mbox{\rm div}\,D{\cal G}(\Psi^*(\theta))-\frac{d}{d\theta}\ln(||\dot{\Psi}^*(\theta)||^2).
\label{A0thetadef}
\end{equation}
Define
\[
\beta=\frac{1}{T}\int_0^T\,A(0,\theta)d\theta=\frac{1}{T}\int_0^T\,\mbox{\rm div}\,D{\cal G}(\Psi^*(\theta))d\theta,
\]
then $\beta <0$ (resp. $\beta >0$) if the limit cycle $\Psi^*(t)$ is linearly stable (resp. unstable).
The periodic change of variables
\begin{equation}
\rho=\sqrt{\varepsilon}\,x\,e^{-\beta\theta}\,e^{\int_0^\theta\,A(0,s)ds}
\label{cov2}
\end{equation}
transforms (\ref{pertavg1}) into 
\begin{equation}
\begin{array}{rcl}
\dot{x}&=&\varepsilon\beta x+\varepsilon^{3/2}V(x,\theta,\varphi,\varepsilon)\\&&\\
\dot{\theta}&=&\varepsilon+\varepsilon^{3/2}W(x,\theta,\varphi,\varepsilon)\\&&\\
\dot{\varphi}&=&1,
\end{array}
\label{pertavg2}
\end{equation}
where all functions are smooth, $T$-periodic in $\theta$ and $2\pi$-periodic in $\varphi$.  System (\ref{pertavg2}) is in the proper form for application of 
Theorem \ref{Ha80:thm23} in the case $n=1$ and $m=2$.  We conclude that there is an invariant two-torus for (\ref{pertavg2}) which is the graph
\[
x=\sigma_{\Psi^*}(\theta,\varphi,\varepsilon),
\]
where $\sigma_{\Psi^*}(\theta,\varphi,\varepsilon)\rightarrow 0$ as $\varepsilon\rightarrow 0$.  From (\ref{cov1}) and (\ref{cov2}), we conclude that (\ref{avpertrot}) has an invariant two-torus
\[
\Psi=\Psi^*(\theta)+\sqrt{\varepsilon}\,\Sigma_{\Psi^*}(\theta,\varphi,\varepsilon),
\]
where
\[
\Sigma_{\Psi^*}(\theta,\varphi,\varepsilon)=J\,{\cal G}(\Psi^*(\theta))
e^{-\beta\theta}\,e^{\int_0^\theta\,A(0,s)ds}\,\sigma_{\Psi^*}(\theta,\varphi,\varepsilon).
\]
We use an argument similar to that used above (for equilibria) to get the conjugate invariant two-tori, and the symmetry property (\ref{stsp2}).
\hfill\qed


\Section{Proof of Theorem \ref{thmtravwave}}

It will be convenient to rewrite (\ref{cmeqs1}) with $\omega=0$ as
\begin{equation}
\begin{array}{rcl}
\dot{\Psi}&=&R_{\varphi}V+\varepsilon\,G^\Psi(\Psi,\varphi,\varepsilon)\\&&\\
\dot{\varphi}&=&\varepsilon\,F^{\varphi}(\Psi,\varphi,\varepsilon)
\end{array}
\label{cmeqs3}
\end{equation}
where $G^\Psi=R_{\varphi}F^\Psi$.
Then for $\varepsilon>0$, the change of variables 
\[
\varphi=\varphi^*+\sqrt{\varepsilon}\,\hat{\varphi}
\]
transforms (\ref{cmeqs3}) into the following (upon dropping the hats):
\[
\begin{array}{rcl}
\dot{\Psi}&=&R_{\varphi^*+\sqrt{\varepsilon}\,\varphi} V+\varepsilon\,G^{\Psi}(\Psi,\varphi^*+\sqrt{\varepsilon}\,\varphi,\varepsilon)\\&&\\
\dot{\varphi}&=&\sqrt{\varepsilon}\,F^{\varphi}(\Psi,\varphi^*+\sqrt{\varepsilon}\,\varphi,\varepsilon).
\end{array}
\]
A Taylor expansion of these equations yields
\begin{equation}
\begin{array}{rcl}
\dot{\Psi}&=&R_{\varphi^*}V+\sqrt{\varepsilon}R'_{\varphi^*}V\varphi+\varepsilon [G^{\Psi}(\Psi,\varphi^*,0)+\frac{1}{2}R''_{\varphi^*}V\varphi^2]+\varepsilon^{\frac{3}{2}} Q_1(\Psi,\varphi,\varepsilon)\\&&\\
\dot{\varphi}&=&\sqrt{\varepsilon}F^{\varphi}(\Psi,\varphi^*,0)+\varepsilon F^{\varphi}_{\varphi}(\Psi,\varphi^*,0)\varphi+\varepsilon^{\frac{3}{2}}[F^{\varphi}_{\varepsilon}(\Psi,\varphi^*,0)+\frac{1}{2}F^{\varphi}_{\varphi\varphi}(\Psi,\varphi^*,0)\varphi^2]+\\[0.1in]
&&\varepsilon^2 S_1(\Psi,\varphi,\varepsilon)
\end{array}
\label{cmeqs2}
\end{equation}

The Diophantine condition (\ref{diophantine}) ensures that for any smooth function ${\cal H}:\Ttwo\longrightarrow\mathbb{R}$, if we denote by
\[
\langle\,{\cal H}\,\rangle=\frac{1}{4\pi^2}\iint_{\Ttwo}\,H(\Psi)\,d\Psi,
\]
then there exists a smooth solution ${\cal F}(\Psi)$ to the equation
\[
D_{\Psi}{\cal F}(\Psi)R_{\varphi^*}V={\cal H}(\Psi)-\langle\,{\cal H}\,\rangle,
\]
with $\langle\,{\cal F}\,\rangle=0$.  

Recall that $\langle\,F^{\varphi}(\cdot,\varphi^*,0)\,\rangle=M(\varphi^*)=0$.
Consequently, if ${\cal Y}(\Psi)$ is such that 
\[
D_{\Psi}{\cal Y}(\Psi)R_{\varphi^*}V=F^{\varphi}(\Psi,\varphi^*,0),
\]
then the near-identity change of variables
\[
\varphi=\hat{\varphi}+\sqrt{\varepsilon}\,{\cal Y}(\hat{\Psi}),\,\,\,\,\,\,\,\,\,\,\,\,\Psi=\hat{\Psi}
\]
transforms (\ref{cmeqs2}) into the following (upon dropping the hats)
\begin{equation}
\begin{array}{rcl}
\dot{\Psi}&=&R_{\varphi^*}V+\sqrt{\varepsilon}R'_{\varphi^*}V\varphi+\varepsilon [G^{\Psi}(\Psi,\varphi^*,0)+\frac{1}{2}R''_{\varphi^*}V\varphi^2+R'_{\varphi^*}V{\cal Y}(\Psi)]+\varepsilon^{\frac{3}{2}} Q_2(\Psi,\varphi,\varepsilon)\\&&\\
\dot{\varphi}&=&\varepsilon [F^{\varphi}_{\varphi}(\Psi,\varphi^*,0)-D_{\Psi}{\cal Y}(\Psi)R'_{\varphi^*}V]\varphi+\varepsilon^{\frac{3}{2}}[F^{\varphi}_{\varphi}(\Psi,\varphi^*,0){\cal Y}(\Psi)+F^{\varphi}_{\varepsilon}(\Psi,\varphi^*,0)+\\[0.1in]
&&
\frac{1}{2}F^{\varphi}_{\varphi\varphi}(\Psi,\varphi^*,0)\varphi^2-D_{\Psi}{\cal Y}(\Psi)(R'_{\varphi^*}V{\cal Y}(\Psi)+G^{\Psi}(\Psi,\varphi^*,0)+\frac{1}{2}R''_{\varphi^*}V\varphi^2)]+\varepsilon^2 S_2(\Psi,\varphi,\varepsilon).
\end{array}
\label{cmeqs4}
\end{equation}
If we rescale $\varphi\rightarrow\sqrt{\varepsilon}\varphi$, (\ref{cmeqs4}) simplifies as
\begin{equation}
\begin{array}{rcl}
\dot{\Psi}&=&R_{\varphi^*}V+{\varepsilon}R'_{\varphi^*}V\varphi+\varepsilon [G^{\Psi}(\Psi,\varphi^*,0)+R'_{\varphi^*}V{\cal Y}(\Psi)]+\varepsilon^{\frac{3}{2}} Q_3(\Psi,\varphi,\varepsilon)\\&&\\
\dot{\varphi}&=&\varepsilon [F^{\varphi}_{\varphi}(\Psi,\varphi^*,0)-D_{\Psi}{\cal Y}(\Psi)R'_{\varphi^*}V]\varphi+\varepsilon [F^{\varphi}_{\varphi}(\Psi,\varphi^*,0){\cal Y}(\Psi)+F^{\varphi}_{\varepsilon}(\Psi,\varphi^*,0)+\\[0.1in]
&&
-D_{\Psi}{\cal Y}(\Psi)(R'_{\varphi^*}V{\cal Y}(\Psi)+G^{\Psi}(\Psi,\varphi^*,0))]+\varepsilon^{\frac{3}{2}} S_3(\Psi,\varphi,\varepsilon).
\end{array}
\label{cmeqs5}
\end{equation}

A final near-identity change of variables
\[
\varphi=\hat{\varphi}+\varepsilon [ {\cal Z}_1(\hat{\Psi})+\hat{\varphi}{\cal Z}_2(\hat{\Psi}) ],\,\,\,\,\,\,\,\,\,\,\,\,\,\,
\Psi=\hat{\Psi}+\varepsilon {\cal Z}_3(\hat{\Psi})
\]
transforms (\ref{cmeqs5}) into the following (upon dropping the hats)
\begin{equation}
\begin{array}{rcl}
\dot{\Psi}&=&R_{\varphi^*}V+\varepsilon [R'_{\varphi^*}V\varphi+\kappa]+\varepsilon^{\frac{3}{2}} Q_4(\Psi,\varphi,\varepsilon)\\&&\\
\dot{\varphi}&=&\varepsilon [\mu\varphi+\lambda]+\varepsilon^{\frac{3}{2}} S_4(\Psi,\varphi,\varepsilon)
\end{array}
\label{cmeqs6}
\end{equation}
where we recall that $\mu=M'(\varphi^*)\neq 0$, and we define
\[
\kappa=\langle\,G^{\Psi}(\cdot,\varphi^*,0)\,\rangle
\]
and
\[
\lambda=\langle\,
F^{\varphi}_{\varphi}(\cdot,\varphi^*,0){\cal Y}(\cdot)+F^{\varphi}_{\varepsilon}(\cdot,\varphi^*,0)
-D_{\Psi}{\cal Y}(\cdot)(R'_{\varphi^*}V{\cal Y}(\cdot)+G^{\Psi}(\cdot,\varphi^*,0))\,\rangle .
\]
Translation of the variable $\varphi$ gives the final system
\begin{equation}
\begin{array}{rcl}
\dot{\Psi}&=&R_{\varphi^*}V+\varepsilon [R'_{\varphi^*}V(\varphi-\lambda/\mu)+\kappa]+\varepsilon^{\frac{3}{2}} Q(\Psi,\varphi,\varepsilon)\\&&\\
\dot{\varphi}&=&\varepsilon \mu \varphi +\varepsilon^{\frac{3}{2}} S(\Psi,\varphi,\varepsilon)
\end{array}
\label{cmeqs7}
\end{equation}
System (\ref{cmeqs7}) satisfies all the hypotheses of Theorem \ref{Ha80:thm23} above (i.e. Theorem 2.3, \S VII.2 of \cite{Ha80}), and so we conclude that (\ref{cmeqs7}) admits an invariant two-torus
$
\varphi={\cal F}(\Psi,\varepsilon),
$
(asymptotically stable (resp. unstable) if $\mu<0$ (resp. $\mu > 0$))
which corresponds to a hyperbolic invariant two-torus
\[
\varphi={\cal T}(\Psi,\varepsilon),
\]
(with ${\cal T}(\Psi,\varepsilon)\rightarrow\varphi^*$ as $\varepsilon\rightarrow 0$)
for (\ref{cmeqs1}).
This concludes the proof of the theorem.
\hfill
\qed

\end{document}